\documentclass[draftclsnofoot]{IEEEtran}

\usepackage{cite}

\ifCLASSINFOpdf
   \usepackage[pdftex]{graphicx}
   \DeclareGraphicsExtensions{.pdf,.jpeg,.png}
\else
   \usepackage[dvips]{graphicx}
   \DeclareGraphicsExtensions{.eps}
\fi

\usepackage[cmex10]{amsmath}
\usepackage{url}

\begin{document}
\title{Constructing circuit codes \\ by permuting initial sequences}
\author{Ed~Wynn
\thanks{E-mail ed.wynn@zoho.com}%
\thanks{Manuscript dated 8 January 2012.}
}

\markboth{Submitted to arXiv}%
{Wynn: Permuted circuit codes}

\maketitle

\begin{abstract}
Two new constructions are presented for coils and snakes in the hypercube. Improvements
are made on the best known results for snake-in-the-box coils of dimensions 9, 10 and
11, and for some other circuit codes of dimensions between 8 and 13. In the first construction,
circuit codes are generated from permuted copies of an initial transition sequence; the
multiple copies constrain the search, so that long codes can be found relatively
efficiently. In the second construction, two lower-dimensional paths are joined
together with only one or two changes in the highest dimension; this requires a search
for a permutation of the second sequence to fit around the first. It is possible to investigate 
sequences of vertices of the hypercube, including
circuit codes, by connecting the corresponding vertices in an extended
graph related to the hypercube. As an example of this, invertible circuit codes are briefly
discussed.
\end{abstract}


\IEEEpeerreviewmaketitle

\section{Introduction}
\label{Sec:Introduction}
\IEEEPARstart{L}{et}
 $I^d$ be the graph of the $d$-dimensional hypercube.
That is, the vertex set of $I^d$ is $\{0,1\}^d$, and two
vertices are connected by an edge if and only if they differ in exactly one coordinate.
A $d$-dimensional \emph{circuit code} of length $N$ and spread $k\geq1$ is 
a simple circuit $(x_0, x_1, \ldots, x_{N-1}, x_0)$ in $I^d$ with the property that
$D(x_i, x_j) \geq \min(k, j-i, N+i-j)$ for all $0 \leq i < j \leq N-1$
where $D(x_i, x_j)$ is the number of coordinates of $I^d$ in which $x_i$ and $x_j$ differ. 
In other words, if two nodes are part of the circuit and the distance between them is $i \leq k$, 
then they must be connected directly by $i$ transitions in the sequence.

A $d$-dimensional circuit code of spread 2 is here called a \emph{$d$-coil}, following the
terminology of \cite{OEIS_snakes}. Discovering long $d$-coils is known as the snake-in-the-box
problem \cite{Kochut1996}, and $d$-coils have been called snakes \cite{Klee1970}. 
In this work, a \emph{$d$-snake} is a simple path in $I^d$ with spread 2. 

A circuit code can be described by its \emph{transition sequence} $(c_0, c_1, \ldots, c_{N-1})$
, where $c_i$ specifies the coordinate that changes from $x_i$ to $x_{i+1}$ 
(with wraparound modulo $N$).
This paper presents two new constructions of circuit codes and some results of these
constructions.

\section{Permuted circuit codes}
\label{Sec:permutes}

From an \emph{initial sequence} $\mathbf{c}^{(0)}=(c_0^{(0)}, c_1^{(0)}, \ldots, c_{L-1}^{(0)})$,
we define \emph{permuted sequences} $\mathbf{c}^{(p)}$ by $c_i^{(p)}=\pi(c_i^{(p-1)})$ 
for $p \geq 1$ and $0 \leq i \leq L-1$, where $\pi$ is a permutation of
$\{0,1,\ldots,{d-1}\}$.  
A \emph{permuted circuit code} of period $P$ is then constructed as a circuit code
whose transition sequence is an initial sequence followed by $P-1$ permuted sequences: 
$(\mathbf{c}^{(0)}, \mathbf{c}^{(1)}, \ldots, \mathbf{c}^{(P-1)})$.
It is convenient to divide the vertices into the corresponding sequences of length $L$: we define
$x_i^{(p)}$ to be $x_{p L+i}$ 
for $p \geq 0$ and $0 \leq i \leq L-1$.

An example of a permuted circuit code is one of the four longest 6-coils, with $N=26$.  This has $P=2$. 
The initial sequence is $(0,1,2,0,3,1,4,0,2,5,3,1,2)$, and the remainder is a simple permutation (merely swapping 2 and 4): 
$(0,1,4,0,3,1,2,0,4,5,3,1,4)$.  
(The other three 6-coils with $N=26$ can be described using terminology of later sections:
one is asymmetric; one is natural; and one is invertible.)

When the initial vertex $x_0$ is assumed without loss of generality to be $\mathbf{0}$, 
then the \emph{initial leap} is defined to be the coordinate vector of $x_0^{(1)}$.
This vector in $I^d$ will involve a change or no change in each coordinate, 
according to whether the initial sequence has an even or odd number of changes in that coordinate.
So, for the example of an initial sequence of in the previous paragraph, the pair of changes in coordinate 3
cancels out but all other coordinates have odd numbers of changes, so the initial leap is (1,1,1,0,1,1).

An algorithm for constructing permuted circuit codes is as follows:

\begin{itemize}

\item All permutations $\pi$ are generated, up to conjugacy.
\item For each permutation, each of the $2^d$ possible vectors is proposed in turn as a possible 
initial leap.  
The initial leap is proposed before the initial sequence, or even the length of the initial 
sequence, is known.
\item When a permutation and an initial leap have been proposed, then successive vertices 
$x_0^{(1)}$, $x_0^{(2)}$, $\ldots$ can then be deduced: coordinate $i$ of one leap vector is equal 
to coordinate $\pi(i)$ of the next.  These vertices are generated until either the $k$-spread 
condition is violated (in which case the initial leap is rejected)
or the initial vertex is revisited (so that $x_0^{(P)} = x_0^{(0)}$ for some $P$).  In this way, 
the period $P$
of a permuted circuit code can be deduced from its permutation and its initial leap.
To return to the example above, a permutation 01(24)35 and an initial leap $(1,1,1,0,1,1)$
would be considered.  The permuted leap is the same as the initial leap, which returns to the
initial vertex with $P=2$.
\item The \emph{skeleton} of initial vertices $(x_0^{(0)},  x_0^{(1)}, \ldots ,  x_0^{(P-1)},  x_0^{(0)})$
is compared to previous skeletons from the same permutation, and duplicates (up to isomorphism 
in $I^d$) are rejected.
Effectively, this is a test whether the proposed initial leap vector is equal to 
$\sigma(\mathbf{v})$, where  $\mathbf{v}$ is a leap vector that has already been proposed, 
and $\sigma$ is a permutation that commutes with $\pi$.  One method for conducting this test is 
described in 
Section \ref{Sec:graphs}.

\item For a suitable permutation and initial leap, an exhaustive search with backtracking is then 
conducted for initial sequences that link $x_0^{(0)}$ to $x_0^{(1)}$.  
Whenever a new change in $c_i^{(0)}$ is proposed, 
a new vertex  $x_{i+1}^{(0)}$ can be tested against the spread-$k$ condition.
Also, the equivalent changes and equivalent vertices in the permuted sequences can be deduced and tested.
Therefore, the backtracking search is subject to many constraints.  
If an initial sequence is found that reaches $x_0^{(1)}$, then each permuted sequence also reaches
the start of the next sequence, and this finally defines $L$ and the whole coil of $N$ vertices, 
with $N=LP$.

\end{itemize}

To summarise this algorithm: an initial leap defines the leap from the starting-point of the initial
sequence, $x_0^{(0)}$, to the starting-point of the first permuted sequence, $x_0^{(1)}$.  
A permutation defines how the initial leap changes into subsequent leaps.  
The algorithm looks for successful combinations of initial leap and permutation, 
which define a skeleton of $P$ starting-points that lead back to the initial vertex.
Then a search is made for a transition sequence from $x_0^{(0)}$ to $x_0^{(1)}$.  During this search,
 each proposed transition defines a new change in each permuted sequence.  
There are two main advantages of this algorithm over a simple backtracking search: permuted copies of every proposed transition are additional constraints; 
and the initial sequence can be short but still produce a long coil.

Another example is supplied, again in $d=6$.  Consider a permutation (345201) and an initial leap of
(0,1,1,0,0,0).  Five permuted leaps are then (0,0,0,0,1,1), (1,1,0,0,0,0), (0,0,0,1,1,0), (1,0,1,0,0,0)
and (0,0,0,1,0,1).  It can be seen that the effect of these six leaps is to return to the original vertex,
so this defines a skeleton with $P=6$.  It turns out that a transition sequence (0,1,2,0), which
accomplishes the initial leap, is compatible with the permuted sequences (3,4,5,3), (2,0,1,2) and so on
that accomplish the permuted leaps, and a permuted coil of $N=24$ is formed: 
(0,1,2,0,3,4,5,3,2,0,1,2,5,3,4,5,1,2,0,1,4,5,3,4).

In the first step, permutations of $\{0,1,\ldots,d-1\}$ are considered conjugate
if they have the same set of cycle lengths.  Therefore, a method of generating 
non-conjugate permutations is to consider all partitions of the integer $d$.  
For each partition, a cycle is generated with that length.
To generate the partitions, Algorithm 7.2.1.4P in \cite{KnuthVol4A_standalone} would be suitable.
This generates partitions in reverse lexicographic order, from `$d$' to `$11\ldots1$'.
It was sometimes found to be efficient to concentrate on permutations with long cycles, when
an exhaustive search was prohibitive.  Therefore a variant of this algorithm was developed,
to generate partitions in lexicographic order.

\section{Examples of permuted circuit codes}
\label{Sec:results}

For $d=10$ and 11, the construction algorithm in the previous section has been used to produce 
$d$-coils of lengths 348 and 640, longer than the previously known longest, 344 and 630 
\cite{CasellaPotter2005}.
For $d$-dimensional circuit codes of spread 3, with $d=10$ and 11,
the construction produces lengths of 100 and 160, longer than the previous known longest, 86 
\cite{PatersonTuliani1997} and 154 \cite{ZinovikKroeningChebiryak2008}.
These new circuit codes are detailed in the Appendix.

\medskip
For all $d \geq 2$, the transition sequence $(0, 1, \ldots,\linebreak[1] {d-1},\linebreak[1] 0, 1, \ldots,\linebreak[1] {d-1})$
defines a $d$-dimensional circuit code of length $2d$ and spread $d$.  This can be
described as a permuted circuit code with a length-1 initial sequence, $\mathbf{c}^{(0)} = (0)$,
with $\pi:i \mapsto (i+1) \mod d$.

\medskip
For all $d \geq 3$, a permuted $d$-coil $C$ of length $2d$ is defined by $\mathbf{c}^{(0)} = (1,0)$
and $\pi:i \mapsto (i+1) \mod d$.  This coil
contains all $d$ neighbors of a vertex.  (The vertex itself is not part of the coil, of course.)
This is most easily seen by taking the first vertex of the circuit to be $\{1,0,0,\ldots,0\}$, with a single
1-coordinate in coordinate 0.  Each pair of changes:

\indent 1 0 \noindent

\indent 2 1 \noindent

\indent 3 2 \noindent

\indent etc. \noindent

\noindent adds a new 1-coordinate and cancels the previous one, so all vertices with a single 
1-coordinate are in the coil $C$.  These are precisely the $d$ neighbors of $\{0, \ldots, 0\}$.  
Any vertex linked in $C$ to one of these neighbors has two 1-coordinates, and must therefore be 
linked to the neighbor with the other 1-coordinate to preserve the spread of 2.  Thus this is the 
unique $d$-coil, up to isomorphism in $I^d$, where a vertex not in the coil has all $d$ neighbors 
in $C$.  For $d=3$, this $d$-coil is identical to the one in the previous paragraph.

\medskip
A special case of a permuted circuit code is where the initial sequence is repeated
once, unchanged: $(c_0, \ldots, c_{L-1},\linebreak[1] c_0, \ldots, c_{L-1})$.  With this property,
a $d$-coil may be called `natural' \cite{Kautz1958} or `symmetric' \cite{Singleton1966}.  
These coils can be regarded as permuted circuit codes with period 2 and the identity permutation.
An exhaustive search has been made for 8-coils of this type, again taking the approach
of proposing trial vectors for the initial leap to $x_L$.  This puts additional constraints 
on the search that starts at $x_0$.  The longest results have length 94 (compared to 96 for
the longest known general 8-coil \cite{PatersonTuliani1997}); an example is given in the Appendix.

\medskip
The definition of a permuted circuit code can be extended to allow the final permuted
sequence to be truncated; the transition sequence is then \[(\mathbf{c}^{(0)}, \mathbf{c}^{(1)}, 
\ldots,\linebreak[1] \mathbf{c}^{(P-2)},\linebreak[1] c_0^{(P-1)}, c_1^{(P-1)},\ldots,c_i^{(P-1)})\] with $i<L-1$.
An example is $d=5$, $\mathbf{c}^{(0)} = (1,0,3)$ and $\pi:i \mapsto (i+1) \mod d$, which 
ends with $N=10$.  When the algorithm was modified to find examples of these special cases,
the results were generally shorter than those from full permuted repetitions with similar
leap periods.

\section{Computation times}
\label{Sec:timings}

Given a permutation and an initial leap, the search for an initial sequence may be 
highly constrained, because every change $c_i^{(0)}$ defines other changes $c_i^{(1)}$ etc.,
and all the new occupied vertices must avoid all other vertices in the coil, with spread $k$.
Also, if $P$ is large, then only a short initial sequence is needed to produce a long coil.
These considerations can make the searches relatively quick.

Example computation times are stated as CPU time for a single processor on an Intel
Q8200 Core2 Quad 2.33GHz processor, running the gcc 3.4.5 compiler in Microsoft
Vista.

An exhaustive search for permuted 9-coils with leap periods $P \geq 12$ took 1 minute; the
longest result has length 180 with $P=12$. An exhaustive search for permuted 10-coils
with leap periods $P>12$ took 80 minutes; the longest result has length 320 with $P=16$.
The successful search with $P=12$ took 4120 minutes (2.8 days). An exhaustive search
for permuted 11-coils with leap periods $P>22$ took 1320 minutes (0.9 days); the
longest result has length 576 with $P=24$. An exhaustive search for 11-coils with $P=22$ took 6.9
weeks, but produced its first length-638 result after approximately 2 days.  The length-640
11-coil mentioned in Section \ref{Sec:permutes} was found by a restricted search of
$P=20$.

The computational times taken by searches for initial sequences can be compared to
the time taken by an exhaustive search for the longest 7-coils, resulting in length 48,
using the method of Section \ref{Sec:direct}.
This search took 3480 minutes (2.4 days) on the computer mentioned above. 
This time is quoted because of the difficulty of like-for-like comparisons with the time
taken by the first reported exhaustive search \cite{Kochut1996}. An exhaustive search for the 
longest 8-coil is clearly prohibitive using current methods.

The search for permuted $d$-coils becomes increasingly difficult for $d$ larger than 12:
large leap periods are relatively quickly found to be unsuccessful, and short leap
periods take prohibitively long times to search exhaustively. This construction is not
expected to have a large useful range of $d$ without modification.

\section{Use of graphs to compare and investigate sequences of hypercube vertices}
\label{Sec:graphs}

In Section \ref{Sec:permutes}, it was efficient to test the skeleton of known vertices 
$(x_0^{(0)}, x_0^{(1)}, \ldots, x_0^{(P-1)}, x_0^{(0)})$, which
represented the starts of all the initial and permuted sequences.  The test was whether this
skeleton was isomorphic with any previously considered skeleton.

In general, it is useful to be able to compare sequences of vertices.  For example, in searching
for coils, it is efficient to
reject initial sequences that are isomorphic to others that have already been tried. 
In the current work, comparisons such
 as these have been implemented by representing the vertex sequences in a graph
related to the hypercube.

\begin{figure}[!t]
\centering
\includegraphics[width=2.5in]{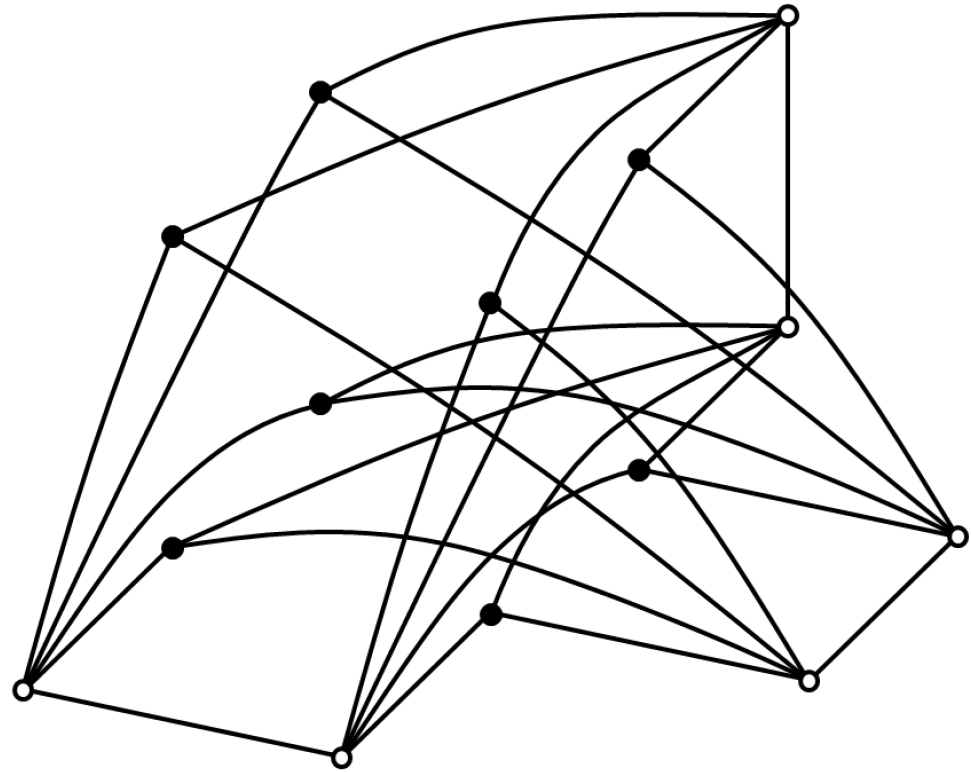}
\caption{The extended graph representing the cube $I^3$.  Filled circles represent original
vertices (in the positions of a plane projection of a geometric cube);
 open circles represent coordinate vertices.}
\label{fig:graph}
\end{figure}

When the hypercube is regarded as a graph, the relationship between the vertices
are defined by connecting each vertex to its neighbours.  Therefore, a sequence of vertices
cannot generally then be defined by simple connections in the resulting graph,
because all permissible connections have already been made.
Therefore, an extended
graph is used here, where the vertices representing those of the hypercube  are not
connected to each other, except when they are linked in a sequence. 
Additional coordinate vertices are used to define the relationship between the original verticles.
The $d$-dimensional extended graph contains $2^d$ original vertices and $2d$ coordinate
vertices. The coordinate vertices come in pairs, each pair representing the 0- and
1-coordinates in one of the $d$ dimensions. Each original vertex is connected to $d$
coordinate vertices, one from each pair, with the choice of 0 or 1 defined by the
relevant coordinate of the corresponding vertex in the hypercube.
 The vertices in each coordinate pair are connected to each other
in the extended graph. The extended graph for $d=3$ is shown in Figure \ref{fig:graph}.

Every automorphism of the hypercube has an equivalent automorphism in
the equivalent extended graph, with original vertices mapping only onto original
vertices, and coordinate vertices onto coordinate vertices.  For example, a reflection 
in one coordinate of the
hypercube is equivalent to swapping the corresponding pair of coordinate vertices. A
sequence of vertices in the hypercube can be represented in the extended graph by
connecting original vertices; sequences that are equivalent (up to symmetries of the
hypercube) are isomorphic in the extended graph. Instances of the extended graph
with connected sequences can be compared efficiently in programs using, for
example, the NAUTY software \cite{McKay1981}: two extended graphs are built up, and NAUTY
then converts them into canonical forms. If and only if these canonical forms are
identical, then the two vertex sequences in the hypercube are equivalent.

There are `brute-force' approaches for comparing vertex sequences in the hypercube.
For example, simple paths can be compared by renumbering the coordinates so that the
transition sequence has its lowest possible lexicographical order, if necessary evaluating
the forwards and backwards versions and selecting the lower of the two. 
 Simple circuits can be compared by similarly
renumbering the coordinates for each starting point, in both directions, and selecting the
lowest of all these renumbered sequences. The extended graph method is generally slower
than these brute-force approaches for simple paths and circuits, but it is more flexible:
for example, it can be applied to chains or skeletons of vertices that are not connected in the 
hypercube.

The extended graph can also be used to study symmetries of paths or circuits; graph
analysis software such as NAUTY can report the orbits of the vertices of an extended
graph with connected original vertices. While the term `symmetric' has sometimes been used for
period-2 permuted $d$-coils with the identity permutation (as mentioned above), this is
not recommended, because there are many other possible symmetries. For example,
any permuted sequence will have order-$P$ symmetry.  The description \emph{natural} is
therefore preferable, although it is not entirely clear how it was derived.

As an example of symmetry in hypercube paths, a $d$-coil can be defined as \emph{invertible}
if the two oriented circuits are isomorphic to each other. 
In the isomorphism, one pair of vertices or changes must remain fixed; we may therefore distinguish
between \emph{vertex-fixed} or \emph{change-fixed} inversions.  In different isomorphisms,
a coil might have both types; an example is the highly symmetrical coil
$(0, 1, \ldots, {d-1},\linebreak[1] 0, 1, \ldots, {d-1})$, but no other examples are known.

Examples of invertible coils are the three longest 5-coils, with $N=14$: two are vertex-fixed invertible and the other is
change-fixed invertible.  One vertex-fixed invertible coil is shown below and repeated in reverse order:

0123142 1023124

4213201 2413210.

\noindent It can be seen that the two presentations are isomorphic via a permutation of coordinates (04)(12)3, 
so the initial vertex is fixed by this inversion.  The same would not be true if the transition sequence
were cycled to start at any other vertex except the `opposite' vertex (at the position of the spaces in the lines above).
So, for example, there is no permutation of coordinates that maps the transition sequence onto its inverse if the initial vertex is moved by one place:

12314210231240

04213201241321.

\noindent The other vertex-fixed invertible 5-coil with $N=14$ is natural:

0123042 0123042

2403210 2403210.

\noindent The change-fixed invertible 5-coil with $N=14$ is shown below and repeated in reverse order, while fixing the two changes:

0 120324 0 123024

0 420321 0 423021.

Using the extended graph technique, we find that all
$d$-coils with $d \leq 5$ are invertible, but longer invertible $d$-coils appear to be rare. 
Of the 4 maximal 6-coils ($N=26$), only 1 is invertible.
The usefulness of the extended graph technique can be illustrated by asking questions:
which of the following two 6-coils is invertible, and is it vertex- or change-fixed?

01203143052351035230420135

01203104201350120310421035.

Of the 758 maximal 7-coils
($N=48$), only 37 are invertible; all except one of these have vertex-fixed inversions.
None of the longer $d$-coils mentioned in this paper are
invertible.

\section{Construction using lower-dimensional snakes}
\label{Sec:other}

An unrelated construction is briefly mentioned here. If the transition sequences
$\mathbf{b}$ and $\mathbf{c}$ define two $d$-snakes, then the transition sequence 
$(\mathbf{b}, d, \mathbf{c})$ may define a $(d+1)$-snake, or
$(\mathbf{b}, d, \mathbf{c}, d)$ may define a $(d+1)$-coil.

Given two $d$-snakes, it will generally be necessary to permute one of the change
sequences to search for a successful combination. An efficient way to do this is to
add vertices of the permuted second snake to the combined sequence one by one,
generating each element of the permutation only when required, and proceeding only
when a suitable element can be found. Algorithm 7.2.1.2X in \cite{KnuthVol4A_standalone} can be used,
because it can be used to generate incomplete permutations exactly as required.
The steps in the method are therefore as follows:
\begin{itemize}
\item Start with a snake consisting of the first sequence.
\item Consider each change in the second sequence, in order.  If it is a change number that has 
not yet been assigned a permutation,
generate a new permutation of that change number.
\item Append the permutation of each change to the current snake.
\item If the new snake disobeys the required spread condition, backtrack to the most recently generated
permutation and generate a new permutation.  If there are no more possibilities, remove the assigned 
permutation and backtrack to the next most recent permutation, and so on.
\end{itemize}

During the exhaustive search for 7-coils, 7-snakes were recorded. Pairs of these were
combined into 8-snakes, which were then similarly combined. For $d = 9$, this
resulted in a coil of length 188 and a snake of length 190, longer than the previously
known longest, 180 and 188 \cite{TuohyPotterCasella2007}. The computation times for attempting to 
combine 9-snakes are not prohibitive, but there is a very large number of candidates of suitable
lengths to form long 10-coils. This construction is not expected to be widely useful for large $d$.

\section{Results from direct searches}
\label{Sec:direct}

Backtracking searches were conducted for circuit codes.  Partial sequences were rejected if any 
subsequence, running forwards or backwards, could be renumbered to a lower number than the 
starting subsequence.  This is not an exhaustive search for snakes, but it avoids the wasted effort 
of finding circuits from multiple starting points.  Every subsequence was tested at every step; 
more efficient strategies may well be available.

Exhaustive searches confirmed the optimality of known sequences \cite{PatersonTuliani1998}:
length 46 for $d = 10$, spread 4; and length 40 for $d = 11$, spread 5.  Additionally, circuit 
codes of lengths 
58, 
  58, and 
    50 were found 
for $d = 9$, spread 3, 
  for $d = 12$, spread 5, and 
    for $d = 13$, spread 6
 --- longer than the previous known longest, 
56 \cite{ZinovikKroeningChebiryak2008},
  56 \cite{ChebiryakKroening2008}, and
    48 \cite{PatersonTuliani1998} respectively.

\section{Conclusions}
\label{Sec:conclusions}

Two methods are presented in this paper for searching for circuit codes.  Both methods attempt
to reduce the combinatorial explosion of the search by adding constraints.  In the first method,
permuted sequences are used; the entire sequence is assembled from multiple permuted copies
of a shorter sequence.  In the second method, the search is effectively for two lower-order
sequences that can be combined with only one pair of changes in the highest coordinate.
It is remarkable that such constrained searches can be competitive with more
general searches, but they have produced new records in several cases that are presented here.

A method is presented for detecting the symmetries of circuit codes through
the use of efficient tools for analysing graphs.  Circuit codes are defined as circuits on the hypercube graph;
the nodes are already linked together before a circuit is specified.  Therefore, it is useful to define 
an extended graph that includes nodes that are closely related to the hypercube's nodes but which are not
initially linked.  Symmetries equivalent to the hypercube's symmetries are brought about through
additional vertices.  Examples of the symmetries that can be detected using this method are two
different inversions.

\appendix[Details of circuit codes]

Some transition sequences are given in full;
others can be deduced by permuting the initial sequences.

\begin{itemize}

\item Coil, spread 2, $d=10$: initial sequence 0189720846 9847685740 278968076, $L=29$,
permutation (123450)(786)9, $P=12$, $N=348$.  The full sequence:

01897208469847685740278968076
12698316579658760851386976187
23796427089706871602467987268
34897538169817682713578968376
45698046279628763824086976487
50796157389736874635167987568
01897208469847685740278968076
12698316579658760851386976187
23796427089706871602467987268
34897538169817682713578968376
45698046279628763824086976487
50796157389736874635167987568.

\item Coil, spread 2, $d=11$: initial sequence 0168763891  305943A671  35127A237, $L=29$,
permutation (123456789A0), $P=22$, $N=638$.

\item Coil, spread 2, $d=11$: initial sequence A04A82A73A 26A38A27A4 8162A648A4 02, $L=32$,
permutation (1234567890)A, $P=20$, $N=640$.  
In this coil, it is noteworthy that more than one quarter of the changes are in a 
single coordinate.  The full sequence:

A04A82A73A26A38A27A48162A648A402
A15A93A84A37A49A38A59273A759A513
A26A04A95A48A50A49A60384A860A624
A37A15A06A59A61A50A71495A971A735
A48A26A17A60A72A61A82506A082A846
A59A37A28A71A83A72A93617A193A957
A60A48A39A82A94A83A04728A204A068
A71A59A40A93A05A94A15839A315A179
A82A60A51A04A16A05A26940A426A280
A93A71A62A15A27A16A37051A537A391
A04A82A73A26A38A27A48162A648A402
A15A93A84A37A49A38A59273A759A513
A26A04A95A48A50A49A60384A860A624
A37A15A06A59A61A50A71495A971A735
A48A26A17A60A72A61A82506A082A846
A59A37A28A71A83A72A93617A193A957
A60A48A39A82A94A83A04728A204A068
A71A59A40A93A05A94A15839A315A179
A82A60A51A04A16A05A26940A426A280
A93A71A62A15A27A16A37051A537A391.

\item Coil, spread 3, $d=10$: initial sequence 26014, $L=5$, permutation (1234567890), $P=20$,
$N=100$.

\item Coil, spread 3, $d=11$: initial sequence 0A184A5234, $L=10$, permutation
(12345670)(98)A, $P=16$, $N=160$.

\item Coil, spread 3, $d=11$: initial sequence 0623184A, $L=8$, permutation (1234567890)A,
$P=20$, $N=160$.

\item Natural coil, $d=8$: transition sequence
0314035046  0340745135  6253157407  5305670517  0317436 twice, $N=94$.

\item Coil, $d=9$, from construction in Section \ref{Sec:other}: transition sequence
0123043254  2134256352  1324532105  1245231524  6142315712  3152413210  4213245321  3461235421  
3253045213  2458032105  1245231524  6142312541  2304325421  3425635213  4732134253  1230523125  
4123156321  4523124105  42312548, $N=188$.

\item Snake, $d=9$, from construction in Section \ref{Sec:other}: transition sequence
0120314021  0541021432  1026431450  4134210431  4501432731  2014301263  2143053102  3053145036  
0431402143  1046806104  3145014310  6302143203  5043203145  3654031405  4375314021  4310451341  
0214316504  5314504120  4501430540, $N=190$.

\item Coil, spread 3, $d=9$: 0123041502 1603570132 4038175014 5671536012 3674563017 60581735, $N=58$.

\item Coil, spread 5, $d=12$: 0123450617 2803196A04 72160548B7 014A836105 82A9167854 0613A84B, $N=58$.

\item Coil, spread 6, $d=13$: 0123456071 82930A142B 9C630529A7 60124A8305 629B4C5A, $N=50$.

\end{itemize}

\section*{Acknowledgment}
\label{Sec:Acknowledgment}

The author is grateful for the encouragement of
Prof. Don Potter at the University of Georgia, USA and Dr. Matthew Fayers at
Queen Mary and Westfield College, University of London, UK.  
Prof. Potter's website \url{www.ai.uga.edu/sib}
is recommended.

\ifCLASSOPTIONcaptionsoff
  \newpage
\fi

\bibliographystyle{IEEEtran}
%

\end{document}